\documentclass[a4paper,11pt]{amsart}
\usepackage{mathrsfs}
\usepackage{bbm}
\usepackage{latexsym, amssymb,amsmath,amsthm}
\usepackage[all, knot]{xy}
\xyoption{arc}

\def \To{\longrightarrow}

\def \gr{\operatorname{gr}}
\def \Hom{\operatorname{Hom}}

\def \Rep{\operatorname{Rep}}
\def \Corep{\operatorname{Corep}}
\def \c{\mathbb{C}}
\def \C{\mathcal{C}}
\def \D{\Delta}
\def \d{\delta}

\def \e{\varepsilon}
\def \M{\mathrm{M}}
\def \N{\mathbb{N}}
\def \S{\mathcal{S}}
\def \Z{\mathbb{Z}}
\def \unit{\epsilon}

\numberwithin{equation}{section}

\newtheorem{theorem}{Theorem}[section]
\newtheorem{lemma}[theorem]{Lemma}
\newtheorem{proposition}[theorem]{Proposition}
\newtheorem{corollary}[theorem]{Corollary}

\newtheorem{convention}[theorem]{Convention}

\begin{document}
\title[QUASI-QUANTUM GROUPS AND TENSOR CATEGORIES]{QUIVERS, QUASI-QUANTUM GROUPS AND \\ FINITE TENSOR CATEGORIES}
\author{Hua-Lin Huang}
\address{School of Mathematics, Shandong University, Jinan, Shandong
250100, China} \email{hualin@sdu.edu.cn}
\author{Gongxiang Liu}
\address{Department of Mathematics, Nanjing University, Nanjing, Jiangsu 210093, China}
\email{gxliu@nju.edu.cn}
\author{Yu Ye}
\address{Department of Mathematics, University of Science and Technology of China,
 Hefei, Anhui 230026, China}
\email{yeyu@ustc.edu.cn}
\date{}
\maketitle

\begin{abstract}
We study finite quasi-quantum groups in their quiver setting
developed recently by the first author. We obtain a classification
of finite-dimensional pointed Majid algebras of finite
corepresentation type, or equivalently a classification of
elementary quasi-Hopf algebras of finite representation type, over
the field of complex numbers. By the Tannaka-Krein duality
principle, this provides a classification of the finite tensor
categories in which every simple object has Frobenius-Perron
dimension 1 and there are finitely many indecomposable objects up to
isomorphism. Some interesting information of these finite tensor
categories is given by making use of the quiver representation
theory.


\vskip 5pt

\noindent{\bf 2000 MR Subject Classification} \ \ 16W35, 18D20,
16G20, 16G60
\end{abstract}

\section{Introduction}

This paper is devoted to the classification of finite quasi-quantum
groups and the associated representation theory, whence the
classification of finite tensor categories \cite{eo}, within the
quiver setting developed recently in \cite{qha1,qha2}. The notion of
quasi-Hopf algebras was introduced by Drinfeld \cite{d2} in
connection with the Knizhnik-Zamolodchikov system of equations from
conformal field theory. The definition of quasi-Hopf algebras is not
selfdual, so there is a dual notion which is called Majid algebra
after Shnider-Sternberg \cite{ss}. In accordance with Drinfeld's
philosophy of quantum groups \cite{d1}, we understand both of these
mutually dual algebraic structures in the framework of quasi-quantum
groups.

We focus on finite-dimensional pointed Majid algebras, or
equivalently elementary quasi-Hopf algebras. Within this
restriction, we can take full advantage of the quiver techniques to
tackle the problems of classification and representation theory.
Recall that by pointed it is meant that the simple subcoalgebras of
the underlying coalgebras are one-dimensional. Dually, by elementary
it is meant that the underlying algebras are finite-dimensional and
their simple modules are one-dimensional.

Quivers are oriented diagrams consisting of vertices and arrows. Due
to Gabriel \cite{gab1,gab2}, in the early 1970's quivers and their
representations became widespread first in the representation theory
of associative algebras. Nowadays these notions show up in various
areas of mathematics and physics. What we rely on is still their
combinatorial behavior and hence very handy applications in the
study of algebraic structures and representation theory. In
connection with Hopf algebras and quantum groups, Hopf quivers
\cite{cr2} and covering quivers \cite{gs} were introduced, see also
\cite{c1, cr1, green, chyz, oz} for related works. For the quiver
setting of the broader class of quasi-quantum groups, it turns out
that there is nothing new other than the Hopf quivers. In principle,
it is shown in \cite{qha1,qha2} that pointed Majid algebras and
elementary quasi-Hopf algebras can be constructed on Hopf quivers
exhaustively with a help of the projective representation theory of
groups and a proper deformation theory.

According to the well-known Tannaka-Krein duality principle (see for
instance \cite{majid,ce}), finite quasi-quantum groups are deeply
related to finite tensor categories. More precisely, the
representation categories of finite-dimensional quasi-Hopf algebras
and the corepresentation categories of finite-dimensional Majid
algebras are finite tensor categories; conversely, finite tensor
categories with some mild conditions are obtained in this way.

In recent years finite tensor categories and finite quasi-quantum
groups have been intensively studied by Etingof, Gelaki, Nikshych,
Ostrik, and many other authors. In \cite{eno}, the fusion
categories, that is, the semisimple finite tensor categories, are
investigated in depth and a number of general properties are
obtained. A systematic study of not necessarily semisimple finite
tensor categories initiated in \cite{eo}, and some classification
results were obtained in \cite{eg1,eg2,g,eg3} through concrete
constructions of elementary quasi-Hopf algebras.

The aim of this paper is to classify finite-dimensional pointed
Majid algebras of finite corepresentation type, or equivalently
elementary quasi-Hopf algebras of finite representation type. The
obvious motivation for this is two-fold. On the one hand, the
algebras of finite representation type are very important in the
representation theory of associative algebras. The study of such
algebras has been a central theme in the area all along. Given an
interesting class of algebras, one is always tempted to classify
those of finite representation type. On the other hand, their
associated representation categories are the finite tensor
categories in which there are only finitely many non-isomorphic
indecomposable objects. Such finite tensor categories are simplest
after the semisimple ones. It is natural to pay prior attention to
these finite tensor categories with very good finiteness property.


The paper is organized as follows. In Section 2 we recall the quiver
setting of quasi-quantum groups. Section 3 is devoted to the
classification of finite-dimensional pointed Majid algebras of
corepresentation type. In Section 4 we investigate finite tensor
categories by making use of quiver representation theory.

Throughout, we work over the field $\mathbb{C}$ of complex numbers
for simplicity. For the convenience of the exposition, we deal
mainly with pointed Majid algebras and mention briefly the situation
of elementary quasi-Hopf algebras. About general background
knowledge, the reader is referred to \cite{ass, ars} for quivers and
representation theory of algebras, to \cite{kassel, majid, ss} for
quasi-quantum groups, and to \cite{bk, ce} for tensor categories.

\section{Hopf Quivers and Quasi-Quantum Groups}

In this section we recall the quiver framework of quasi-quantum
groups for the convenience of the reader.

\subsection{Hopf Quivers}

A quiver is a quadruple $Q=(Q_0,Q_1,s,t),$ where $Q_0$ is the set of
vertices, $Q_1$ is the set of arrows, and $s,t:\ Q_1 \longrightarrow
Q_0$ are two maps assigning respectively the source and the target
for each arrow. A path of length $l \ge 1$ in the quiver $Q$ is a
finitely ordered sequence of $l$ arrows $a_l \cdots a_1$ such that
$s(a_{i+1})=t(a_i)$ for $1 \le i \le l-1.$ By convention a vertex is
said to be a trivial path of length $0.$

For a quiver $Q,$ the associated path coalgebra $\c Q$ is the
$\c$-space spanned by the set of paths with counit and
comultiplication maps defined by $\e(g)=1, \ \D(g)=g \otimes g$ for
each $g \in Q_0,$ and for each nontrivial path $p=a_n \cdots a_1, \
\e(p)=0,$
\begin{equation}
\D(a_n \cdots a_1)=p \otimes s(a_1) + \sum_{i=1}^{n-1}a_n \cdots
a_{i+1} \otimes a_i \cdots a_1 \nonumber + t(a_n) \otimes p \ .
\end{equation}
The length of paths gives a natural gradation to the path coalgebra.
Let $Q_n$ denote the set of paths of length $n$ in $Q,$ then $\c
Q=\oplus_{n \ge 0} \c Q_n$ and $\D(\c Q_n) \subseteq
\oplus_{n=i+j}\c Q_i \otimes \c Q_j.$ Clearly $\c Q$ is pointed with
the set of group-likes $G(\c Q)=Q_0,$ and has the following
coradical filtration $$ \c Q_0 \subseteq \c Q_0 \oplus \c Q_1
\subseteq \c Q_0 \oplus \c Q_1 \oplus \c Q_2 \subseteq \cdots.$$
Hence $\c Q$ is coradically graded. The path coalgebras can be
presented as cotensor coalgebras, so they are cofree in the category
of pointed coalgebras and enjoy a universal mapping property.

According to \cite{cr2}, a quiver $Q$ is said to be a Hopf quiver if
the corresponding path coalgebra $\c Q$ admits a graded Hopf algebra
structure. Hopf quivers can be determined by ramification data of
groups. Let $G$ be a group and denote its set of conjugacy classes
by $\C.$ A ramification datum $R$ of the group $G$ is a formal sum
$\sum_{C \in \C}R_CC$ of conjugacy classes with coefficients in
$\mathbb{N}=\{0,1,2,\cdots\}.$ The corresponding Hopf quiver
$Q=Q(G,R)$ is defined as follows: the set of vertices $Q_0$ is $G,$
and for each $x \in G$ and $c \in C,$ there are $R_C$ arrows going
from $x$ to $cx.$ It is clear by definition that $Q(G,R)$ is
connected if and only if the the set $\{c \in C | C \in \C \
\text{with} \ R_C \ne 0\}$ generates the group $G.$ For a given Hopf
quiver $Q,$ the set of graded Hopf structures on $\c Q$ is in
one-to-one correspondence with the set of $\c Q_0$-Hopf bimodule
structures on $\c Q_1.$

\subsection{Quasi-Quantum Groups}

We recall explicitly the definitions about Majid algebras only.
Those about quasi-Hopf algebras can be written out in a dual manner.

A dual quasi-bialgebra, or Majid bialgebra for short, is a coalgebra
$(H,\D,\e)$ equipped with a compatible quasi-algebra structure.
Namely, there exist two coalgebra homomorphisms $$\M: H \otimes H
\To H, \ a \otimes b \mapsto ab, \quad \mu: k \To H,\ \lambda
\mapsto \lambda 1_H$$ and a convolution-invertible map $\Phi:
H^{\otimes 3} \To k$ called reassociator, such that for all $a,b,c,d
\in H$ the following equalities hold:
\begin{gather}
a_1(b_1c_1)\Phi(a_2,b_2,c_2)=\Phi(a_1,b_1,c_1)(a_2b_2)c_2,\\
1_Ha=a=a1_H, \\
\Phi(a_1,b_1,c_1d_1)\Phi(a_2b_2,c_2,d_2) \\
 =\Phi(b_1,c_1,d_1)\Phi(a_1,b_2c_2,d_2)\Phi(a_2,b_3,c_3),\nonumber \\
\Phi(a,1_H,b)=\e(a)\e(b).
\end{gather}
Here and below we use the Sweedler sigma notation $\D(a)=a_1 \otimes
a_2$ for the coproduct. $H$ is called a Majid algebra if, moreover,
there exist a coalgebra antimorphism $\S: H \To H$ and two
functionals $\alpha,\beta: H \To k$ such that for all $a \in H,$
\begin{gather}
\S(a_1)\alpha(a_2)a_3=\alpha(a)1_H, \quad
a_1\beta(a_2)\S(a_3)=\beta(a)1_H, \\
\Phi(a_1,\S(a_3), a_5)\beta(a_2)\alpha(a_4)= \\
\Phi^{-1}(\S(a_1),a_3,\S(a_5)) \alpha(a_2)\beta(a_4)=\e(a).
\nonumber
\end{gather}

A Majid algebra $H$ is said to be pointed, if the underlying
coalgebra is pointed. For a given pointed Majid algebra $(H,\D, \e,
\M, \mu, \Phi,\S,\alpha,\beta),$ let $\{H_n\}_{n \ge 0}$ be its
coradical filtration, and $\gr H = H_0 \oplus H_1/H_0 \oplus H_2/H_1
\oplus \cdots$ the corresponding coradically graded coalgebra. Then
$\gr H$ has an induced graded Majid algebra structure. The
corresponding graded reassociator $\gr\Phi$ satisfies
$\gr\Phi(\bar{a},\bar{b},\bar{c})=0$ for all
$\bar{a},\bar{b},\bar{c} \in \gr H$ unless they all lie in $H_0.$
Similar condition holds for $\gr\alpha$ and $\gr\beta.$ In
particular, $H_0$ is a sub Majid algebra and turns out to be the
group algebra $kG$ of the group $G=G(H),$ the set of group-like
elements of $H.$

\subsection{Quiver Setting for Quasi-Quantum Groups}

It is shown in \cite{qha1} that the path coalgebra $\c Q$ admits a
graded Majid algebra structure if and only if the quiver $Q$ is a
Hopf quiver. Moreover, for a given Hopf quiver $Q=Q(G,R),$ if we fix
a Majid algebra structure on $\c Q_0=(\c G,\Phi)$ with
quasi-antipode $(\S,\alpha,\beta),$ then the set of graded Majid
algebra structures on $\c Q$ with $\c Q_0=(\c
G,\Phi,\S,\alpha,\beta)$ is in one-to-one correspondence with the
set of $(\c G,\Phi)$-Majid bimodule structures on $\c Q_1.$
According to \cite{d2}, by transforming the quasi-antipode
$(\S,\alpha,\beta)$ via convolution invertible functionals in
$\Hom_\c (\c G,\c),$ one obtains all the graded Majid algebra
structures on $\c Q$ with $\c Q_0=(\c G,\Phi)$ and an arbitrary
quasi-antipode.

The category of Majid bimodules over a general group with an
arbitrary 3-cocycle is characterized in \cite{qha2} by the
admissible collections of projective representations. Let $G$ be a
group and $\Phi$ a 3-cocycle on $G.$ Denote by $\C$ the set of
conjugacy classes of $G$ and by $Z_C$ the centralizer of one of the
elements, say $g_{_C},$ in the class $C \in \C.$ Let
$\Tilde{\Phi}_C$ be a 2-cocycle on $Z_C$ defined by
\begin{equation}
\Tilde{\Phi}_C(e,f)=\frac{\Phi(e,f,g_{_C})\Phi(ef,f^{-1},e^{-1})\Phi(e,fg_{_C},f^{-1})}
{\Phi(efg_{_C},f^{-1},e^{-1})\Phi(e,f,f^{-1})}\ .
\end{equation}
Then the category of $(\c G,\Phi)$-Majid bimodules is equivalent to
the product of categories $\prod_{C \in \C} (\c
Z_C,\Tilde{\Phi}_C)\mathrm{-rep},$ where $(\c
Z_C,\Tilde{\Phi}_C)\mathrm{-rep}$ is the category of projective
$\Tilde{\Phi}_C$-representations, or equivalently the left module
category of the twisted group algebra $\c ^{\Tilde{\Phi}_C}Z_C$ (see
\cite{karp}).

Thanks to the Gabriel type theorem in \cite{qha1}, for an arbitrary
pointed Majid algebra $H,$ its graded version $\gr H$ can be
realized uniquely as a large sub Majid algebra of some graded Majid
algebra structure on a Hopf quiver. By ``large" it is meant the sub
Majid algebra contains the set of vertices and arrows of the Hopf
quiver. Therefore, in principle all pointed Majid algebras are able
to be constructed on Hopf quivers. The classification project can be
carried out in two steps. The first step is to classify large sub
Majid algebras of those on path coalgebras. This gives a
classification of graded pointed Majid algebras. The second step is
to perform a suitable deformation process to get general pointed
Majid algebras from the graded ones.

\subsection{Multiplication Formula for Quiver Majid Algebras}

In order to construct graded Majid algebras on Hopf quivers, we need
to compute the product of paths. It is shown in \cite{qha1} that the
multiplication formula can be given via quantum shuffle product
\cite{rosso2}.

Suppose that $Q$ is a Hopf quiver with a necessary $\c Q_0$-Majid
bimodule structure on $\c Q_1.$ Let $p \in Q_l$ be a path. An
$n$-thin split of it is a sequence $(p_1, \ \cdots, \ p_n)$ of
vertices and arrows such that the concatenation $p_n \cdots p_1$ is
exactly $p.$ These $n$-thin splits are in one-to-one correspondence
with the $n$-sequences of $(n-l)$ 0's and $l$ 1's. Denote the set of
such sequences by $D_l^n.$ Clearly $|D_l^n|={n \choose l}.$ For
$d=(d_1, \ \cdots, \ d_n) \in D_l^n,$ the corresponding $n$-thin
split is written as $dp=((dp)_1, \ \cdots, \ (dp)_n),$ in which
$(dp)_i$ is a vertex if $d_i=0$ and an arrow if $d_i=1.$ Let
$\alpha=a_m \cdots a_1$ and $\beta=b_n \cdots b_1$ be paths of
length $m$ and $n$ respectively. Let $d \in D_m^{m+n}$ and $\bar{d}
\in D_n^{m+n}$ the complement sequence which is obtained from $d$ by
replacing each 0 by 1 and each 1 by 0. Define an element
$$(\alpha \cdot \beta)_d=[(d\alpha)_{m+n}.(\bar{d}\beta)_{m+n}] \cdots
[(d\alpha)_1.(\bar{d}\beta)_1]$$ in $\c Q_{m+n},$ where
$[(d\alpha)_i.(\bar{d}\beta)_i]$ is understood as the action of $\c
Q_0$-Majid bimodule on $\c Q_1$ and these terms in different
brackets are put together by cotensor product, or equivalently
concatenation. In terms of these notations, the formula of the
product of $\alpha$ and $\beta$ is given as follows:
\begin{equation}
\alpha \cdot \beta=\sum_{d \in D_m^{m+n}}(\alpha \cdot \beta)_d \ .
\end{equation}

We should remark that, for general Majid algebras, the product is
not associative. So the order must be concerned for the product of
more than two terms.

\begin{convention}
For an arbitrary path $p$ and an integer $n \ge 3,$ let
$p^{\stackrel{\rightharpoonup}{n}}$ denote the product
$\stackrel{n-2}{\overbrace{( \cdots (}}p \cdot p) \cdot \cdots)
\cdot p$ calculating from the left side. For consistency, when $n <
3,$ we still use the notation $p^{\stackrel{\rightharpoonup}{n}}$
although there is no risk of associative problem. Similarly we use
the notation $p^{\stackrel{\leftharpoonup}{n}}$ for the product
calculating from the right side.
\end{convention}

\section{Pointed Majid Algebras of Finite Corepresentation Type}

In this section we give an explicit classification of
finite-dimensional graded pointed Majid algebras of finite
corepresentation type. We start by fixing the Hopf quivers on which
such Majid algebras live. Then we calculate all the possible Majid
bimodules for the construction. Finally we provide the
classification by making use of quiver techniques.

\subsection{Determination of Hopf Quivers}

Recall that a finite-dimensional algebra is defined to be of finite
representation type if the number of the isomorphism classes of
indecomposable finite-dimensional modules is finite. A
finite-dimensional coalgebra $C$ is said to be of finite
corepresentation type if the dual algebra $C^*$ is of finite
representation type. Since finite-dimensional coalgebras and
finite-dimensional algebras are dual to each other, we apply the
known results of algebras to the coalgebra setting without
explanation.

It is well-known that the module category of a finite-dimensional
elementary algebra, or the comodule category of a finite-dimensional
pointed coalgebra, can be visualized as the representation category
of the corresponding bound quiver (see \cite{ass,chin}). Hence the
quiver presentation of algebras or coalgebras can provide important
information for their representation or corepresentation type. This
is the starting point of our classification.

\begin{lemma}
Let $C\neq \c$ be a finite-dimensional pointed coalgebra and assume
that its bound quiver $Q$ is a connected Hopf quiver. Then $C$ is of
finite corepresentation type if and only if its bound quiver
$Q=Q(\Z_n,g)$ where $\Z_n=\langle g|g^{n}=\unit \rangle$ for some
integer $n \ge 1.$
\end{lemma}

\begin{proof} ``$\Rightarrow$"  We denote the vertex group of the quiver $Q$ by $G$ and
assume that $Q=Q(G,R)$ for some ramification datum $R=\sum_{C \in
\C} R_CC$ of $G.$ Note that a Hopf quiver is very symmetric, its
shape is completely determined once the arrows with source $\unit,$
the unit element of $G,$ are known. For our Hopf quiver $Q,$ we
claim that the number of arrows starting from $\unit$ is 1.

Assume otherwise there are at least two arrows with source $\unit.$
Then the Hopf quiver $Q$ either contains the Kronecker quiver
$$ \xy (0,0)*{\circ \ }; (30,0)*{\ \circ}; {\ar (1,1)*{};
(29,1)*{}}; {\ar (1,-1)*{}; (29,-1)*{}}
\endxy $$ as a sub quiver, or contains a sub quiver of the form
$$ \xy (0,0)*{g}; (15,0)*{\unit}; (30,0)*{h}; {\ar (14,0)*{};
(1,0)*{}}; {\ar (16,0)*{}; (29,0)*{}}; \endxy $$ for some $g \ne h.$
For the latter case, by the definition of Hopf quivers, in $Q$ there
are also arrows as the following
$$ \cdots \longleftarrow h^{-2}g^2 \To h^{-1}g^2 \longleftarrow
h^{-1}g \To g. $$ Since $G$ is a finite group, there is a positive
integer $N$ such that $h^{-N+1}g^N$ or $h^{-N}g^N=\unit.$ It follows
that the quiver $Q$ contains a sub quiver whose underlying graph is
a cycle and whose paths are of length less than 2. By the Gabriel
type theorem for pointed coalgebras (see e.g. \cite{chz,chin}), the
bound quiver $Q$ of $C$ contains either of the above two quivers as
a sub quiver, hence the corepresentation category of $C$ contains
either of their representation categories as a sub category.
According to Gabriel's famous classification of quivers of finite
representation type \cite{gab1}, both of these two quivers admit
infinitely many finite-dimensional indecomposable representations.

On the other hand, since the quiver $Q$ is assumed to be connected,
so the number of arrows starting from $\unit$ can only be 1. Assume
that $\unit \To g$ is the unique arrow of $Q$ with source $\unit.$
By the definition of Hopf quiver, the element $g$ itself must
constitute a conjugacy class and generate the group $G$ since $Q$ is
connected. Therefore, such a quiver can only be of the following
form
$$ \xy {\ar (0,0)*{\unit}; (30,-10)*{g}}; {\ar (-30,-10)*{g^{n-1}};
(0,0)*{\unit}}; {\ar (30,-10)*{g}; (3,-10)*{ \cdots  \ }}; {\ar
(-3,-10)*{\ \cdots  }; (-30,-10)*{g^{n-1}}}
\endxy $$ where the set of vertices $\{ \unit, g, \cdots \cdots,
g^{n-1} \}$ constitutes a cyclic group $\Z_n$ of order $n$ and
$Q=Q(\Z_n,g).$

``$\Leftarrow$" Assume that $Q=Q(\Z_n,g)$ is the bound quiver of
$C,$ then $C$ is a large sub coalgebra of $\c Q$. Denote $\c
Q(d):=\oplus _{i=0}^{d-1}\c Q_{i},$ the $d$-truncated sub coalgebra.
Since $C$ is finite-dimensional, $C$ is a sub coalgebra of $\c Q(d)$
for some $d$. It is well-known that $\c Q(d)$ is of finite
corepresentation type \cite{ass, ars} and thus so is $C$.
\end{proof}

Assume that $M$ is a finite-dimensional pointed Majid algebra and
the corresponding Hopf quiver is $Q.$ To avoid the trivial case, we
assume that $Q$ contains at least one arrow. This excludes the
situation for $M$ being cosemisimple. According to the quasi-Hopf
analogue of the Cartier-Gabriel decomposition theorem in
\cite{qha1}, we may assume without loss of generality that the
quiver $Q$ is connected. This is equivalent to saying that the
underlying coalgebra of $M$ is connected.

\begin{corollary} Keep the previous assumption. If $M$
is of finite corepresentation type then the set of group-likes
$G(M)=\mathbb{Z}_{n}$ and $Q=Q(\mathbb{Z}_{n},g)$ for some integer
$n\geq 2$.
\end{corollary}
\begin{proof} By Lemma 3.1, $Q=Q(\mathbb{Z}_{n},g)$ for some $n\geq 1$.
We claim the case $n=1$ would not occur. Otherwise, the quiver $Q$
consists of one vertex and one loop. By the multiplication formula
(2.8), it is easy to see that the sub Majid algebra generated by the
loop is the shuffle algebra in one variable (see e.g. \cite{kassel})
which is of course infinite-dimensional. By the Gabriel type
theorem, this must be contained in $\gr M.$ This implies that $\gr
M,$ hence $M,$ is infinite-dimensional. Contradiction.
\end{proof}

\subsection{Twisted Group Algebras}

From now on we let $Z^n$ denote the Hopf quiver $Q(\Z_n, g),$ which
is known as the basic cycle of length $n.$ In order to classify
graded Majid algebras on $Z^n$ we need to classify $(\c \Z_n,
\Phi)$-Majid bimodule structures on $\c Z^n_1$ for an arbitrary
3-cocycle $\Phi$ on $\Z_n,$ which can be reduced to the
classification of one-dimensional modules over some twisted group
algebra $\c^\sigma \Z_n$ by \cite{qha2}.

Firstly we recall the nontrivial 3-cocycles on $\Z_n.$ It is
well-known that $H^3(\Z_n,\c^*) \cong \Z_n,$ so there are $n$
mutually non-cohomologous 3-cocycles. We give a list after Lemmas
3.3 and 3.4 of \cite{g}. Let $\mathbbm{q}$ be a primitive root of
unity of order $n.$ For any integer $i \in \N,$ we denote by $i'$
the remainder of division of $i$ by $n.$ A list of 3-cocycles on
$\Z_n$ are
\begin{equation}
\Phi_s(g^{i},g^{j},g^{k})=\mathbbm{q}^{si(j+k-(j+k)')/n}
\end{equation}
for all $0 \le s \le n-1$ and $0 \le i,j,k \le n-1.$ Obviously,
$\Phi_s$ is trivial (i.e., cohomologous to a 3-coboundary) if and
only if $s=0.$

For a 3-cocycle $\Phi_s,$ we define a 2-cocycle $\sigma_s$ on $\Z_n$
by (2.7) as follows:
$$\sigma_s(g^{i},g^{j}):=\frac{\Phi_s(g^{i},g^{j},g)\Phi_s(g^{i+j},g^{-j},g^{-i})\Phi_s(g^{i},g^{j+1},g^{-j})}
{\Phi_s(g^{i+j+1},g^{-j},g^{-i})\Phi_s(g^{i},g^{j},g^{-j})}$$ for
all $0\leq i,j\leq n-1$. Consider the associated twisted group
algebra $\c^{\sigma_s}\mathbb{Z}_{n}.$ We denote the multiplication
in $\c^{\sigma_s}\mathbb{Z}_{n}$ by $``\ast"$. Thus
$$g\ast g=\sigma_s(g,g)g^{2}=\frac{\Phi_s(g,g,g)\Phi_s(g^{2},g^{n-1},g^{n-1})\Phi_s(g,g^{2},g^{n-1})}
{\Phi_s(g^{3},g^{n-1},g^{n-1})\Phi_s(g,g,g^{n-1})}g^{2}=\mathbbm{q}^{-s}g^{2}.$$
Denote $\stackrel{i}{\overbrace{g \ast \cdots \ast g}}$ by $g^{\ast
i}$, then we have

\begin{lemma} In $\c^{\sigma_s}\mathbb{Z}_{n}$, we have
\begin{equation} g^{i}=\mathbbm{q}^{(i-1)s}g^{\ast i}
\end{equation} for $1\leq i\leq n$.
\end{lemma}
\begin{proof} Induction on $i$,
\begin{eqnarray*}
g \ast g^{i}&=&\mathbbm{q}^{(i-1)s}g^{\ast (i+1)}=\sigma_s(g,g^{i})g^{i+1}\\
&=&\frac{\Phi_s(g,g^{i},g)\Phi_s(g^{i+1},g^{n-i},g^{n-1})\Phi_s(g,g^{i+1},g^{n-i})}
{\Phi_s(g^{i+2},g^{n-i},g^{n-1})\Phi_s(g,g^{i},g^{n-i})}g^{i+1}\\
&=&\mathbbm{q}^{-s}g^{i+1}
\end{eqnarray*}
This implies that $g^{i+1}=\mathbbm{q}^{is}g^{\ast (i+1)}$.
\end{proof}

Now we consider one-dimensional modules of
$\c^{\sigma_s}\mathbb{Z}_{n}.$ Let the one-dimensional space $V=\c
X$ be a $\c^{\sigma_s}\mathbb{Z}_{n}$-module with action given by
\begin{equation}
g\triangleright X = \lambda X
\end{equation}
for some $\lambda\in \c.$ Then we have

\begin{lemma} $\lambda^{n}=\mathbbm{q}^{s}$.\end{lemma}
\begin{proof} Indeed, $\lambda^{n}X=g^{\ast n}\triangleright X=
\mathbbm{q}^{-(n-1)s}g^{n}\triangleright X=\mathbbm{q}^{s}X.$
\end{proof}

So $\lambda=q^{s}$ for some $q$ which is an $n$-th root of
$\mathbbm{q}.$ Thus $q$ is a primitive root of unity of order $n^2$
when $s \ne 0.$ Therefore when $s \ne 0,$ the set of one-dimensional
$\c^{\sigma_s}\mathbb{Z}_{n}$-modules is in one-to-one
correspondence with the set of $n$-th roots of $\mathbbm{q}.$ When
$s=0,$ we have $\c^{\sigma_s}\mathbb{Z}_{n} = \c \mathbb{Z}_n$ and
the the set of one-dimensional $\mathbb{Z}_n$-module is in
one-to-one correspondence with the set of $n$-th roots of unity.

\subsection{Computation of Majid Bimodules}

Here and below, let $X_{i}$ denote the arrow $g^{i-1}
\longrightarrow g^{i}$ of the Hopf quiver $Z^n$ for $1\leq i\leq n.$
For convenience the subscript of $X_i$ is read modulo $n$ in some
circumstances. By Theorem 3.3 of \cite{qha2}, the $(\c \Z_n
,\Phi_s)$-Majid bimodule structures on $\c Z^n_1$ can be obtained by
extending the $\c^{\sigma_s} \Z_n$-module structures on the
one-dimensional space $\c X_{1}.$

Recall that, for an arbitrary group $G$ and a 3-cocycle $\Phi,$  a
$(kG,\Phi)$-Majid bimodule $M$ is simultaneously a $kG$-bicomodule
and a quasi $kG$-bimodule such that the quasi-module structure maps
are $kG$-bicomodule morphisms. Assume that $M=\bigoplus_{g,h \in G}
\ ^gM^h$ is the decomposition into isotypic components, where
$$^gM^h=\{ m \in M \ | \ \d_{_L}(m)=g \otimes m, \ \d_{_R}(m)=m
\otimes h \} \ .$$ Here we use $(\d_{_L},\d_{_R})$ to denote the
bicomodule structure maps. The quasi-actions satisfy the
quasi-associativity, namely
\begin{gather}
e.(f.m)=\frac{\Phi(e,f,g)}{\Phi(e,f,h)}(ef).m,\\
(m.e).f=\frac{\Phi(h,e,f)}{\Phi(g,e,f)}m.(ef),\\
(e.m).f=\frac{\Phi(e,h,f)}{\Phi(e,g,f)}e.(m.f),
\end{gather}
for all $e,f,g,h \in G$ and $m \in \ ^gM^h$. These equalities will
be used freely.

A $(kG,\Phi)$-Majid bimodule can be associated to an admissible
collection of projective modules as follows. We still let $\unit$
denote the unit element, and let $\C$ denote the set of conjugacy
classes of the group $G.$ For each $C \in \C,$ let $Z_C$ denote the
centralizer of one element in $C,$ say $g(C),$ and $\Tilde{\Phi}_C$
the corresponding 2-cocycle $\Tilde{\Phi}_{_{g(C)}}$ on $\Z_C$ as
defined in (2.7), and $M_C = \ ^{g(C)}M^\unit$ the
$\Tilde{\Phi}_C$-representation of $Z_C$ given by \begin{equation} h
\triangleright m = (h.m).h^{-1}, \quad \forall \ h \in Z_C, \ m \in
M_C. \end{equation} Then $(M_C)_{C \in \C}$ is called the
corresponding admissible collections of projective representations
of $M.$ Conversely, given an admissible collection of projective
representations, one can extend it by a twisted version of induced
representation to a Majid bimodule, see \cite{qha2} for detail. This
provides the category equivalence mentioned previously in Subsection
2.3.

Now let's get back to our situation of the cyclic group $\Z_n.$ Let
$q$ be an $n$-th root of $\mathbbm{q}$ and fix an $\c^{\sigma_s}
\Z_n$-action on $\c X_{1}$ by
\begin{equation}
g\triangleright X_{1}=q^s X_{1}.
\end{equation}
We extend this $\c^{\sigma_s} \Z_n$-module to an $(\c \Z_n,
\Phi_s)$-Majid bimodule on $\c Z^n_1.$ The bicomodule structure is
defined according to the quiver structure, namely, for $1 \le i \le
n,$
\begin{equation} \d_{_L} (X_i) = g^i \otimes X_i, \quad \d_{_R} (X_i) = X_i \otimes g^{i-1}. \end{equation}
For the quasi bimodule, there is no harm to assume that
$g.X_{i}=X_{i+1}$ for $1\leq i\leq n-1.$ With this, we have

\begin{lemma} The following equations
\begin{eqnarray}
g.X_{i}&=&X_{i+1} \quad (1\leq i\leq n-1) \\
g.X_{n}&=&\mathbbm{q}^{s}X_{1} \\
X_{i}.g&=&\mathbbm{q}^{-s}q^{-s}X_{i+1} \quad (1\leq i\leq n)
\end{eqnarray}
define a quasi $\c \Z_n$-bimodule on $\c Z^n_1$ and make it a $(\c
\Z_n,\Phi_s)$-Majid bimodule together with the bicomodule structure
defined by (3.9).
\end{lemma}

\begin{proof} Inductively, we have
$$X_{i}=g^{i-1} . X_{1},\;\;\;\;\textrm{for } 1\leq i\leq n .$$
Thus
$$g . X_{n}=g . (g^{n-1} .
X_{1})=\frac{\Phi_s(g,g^{n-1},g)}{\Phi_s(g,g^{n-1},\unit)}g^{n} .
X_{1}=\mathbbm{q}^{s}X_{1}.$$

We proceed to determine the right quasi-action. On the one hand,
$$g^{n-1}\triangleright X_{1}=\mathbbm{q}^{(n-2)s}g^{\ast n-1}\triangleright X_{1}=\mathbbm{q}^{-2s}q^{(n-1)s}X_{1}
=\mathbbm{q}^{-s}q^{-s}X_{1}.$$ On the other hand, by the relation
between $``\triangleright"$ and the quasi-actions $``."$ given by
(3.7), we have
$$g^{n-1}\triangleright X_{1}=(g^{n-1} . X_{1}) . g=X_{n} . g \ .$$
Thus $X_{n} . g=\mathbbm{q}^{-s}q^{-s}X_{1}$.

Now assume that $X_{n-1} . g= cX_{n}$ for some $c \in \c$. Thus
$$X_{n} . g=(g . X_{n-1}) . g=\frac{\Phi_s(g,g^{n-2},g)}{\Phi_s(g,g^{n-1},g)}
g . (X_{n-1} . g)=c\mathbbm{q}^{-s}g . X_{n}=cX_{1}$$ and $X_{n} .
g=\mathbbm{q}^{-s}q^{-s}X_{1}$. Therefore,
$c=\mathbbm{q}^{-s}q^{-s}$. Inductively, assume that $X_{i} .
g=\mathbbm{q}^{-s}q^{-s}X_{i+1}$ for some $i\leq n-1$. Suppose that
$X_{i-1} . g=cX_{i}$ for some $c\in \c$, then
$$\mathbbm{q}^{-s}q^{-s}X_{i+1}=X_{i} . g=(g . X_{i-1}) . g=g . (X_{i-1} . g)
=cg . X_{i}=cX_{i+1}.$$ Thus $c=\mathbbm{q}^{-s}q^{-s}$.

It is straightforward to verify that the quasi-bimodule structure
maps are bicomodule morphisms. This completes the proof.
\end{proof}

\subsection{Graded Majid Algebras on $Z^n$}

In this subsection we calculate the graded Majid algebra on the
quiver $Z^n$ associated to the $(\c \Z_n, \Phi)$-Majid bimodule
given in Lemma 3.5.

Firstly we need to fix some notations. For any $\hbar \in \c$,
define $l_\hbar=1+\hbar+\cdots +\hbar^{l-1}$ and $l!_\hbar=1_\hbar
\cdots l_\hbar$. The Gaussian binomial coefficient is defined by
$\binom{l+m}{l}_\hbar:=\frac{(l+m)!_\hbar}{l!_\hbar m!_\hbar}$.

In the path coalgebra $\c Z^n$, let $p_{i}^{l}$ denote the path
starting from $g^{i}$ with length $l.$ The index $``i"$ is read
modulo $n$ when there is no risk of confusion. We keep the $(\c
\Z_n,\Phi_s)$-Majid bimodule on $\c Z^n_1$ as in Lemma 3.5 and
consider the associated graded Majid algebra as given in
\cite{qha1}. Let $\c Z^n(s,q)$ denote the resulted Majid algebra on
the path coalgebra $\c Z^n$ and $``\cdot"$ denote its
multiplication. Note that the quasi-antipode
$(\mathcal{S},\alpha,\beta)$ of $\c Z^n(s,q)$ satisfies
$\S(g)=g^{-1}, \ \alpha(g)=1, \ \beta(g)=1/\Phi_s(g,g^{-1},g)$ for
all $g \in \Z_n.$

\begin{lemma} For any natural number $l$, we have
$$X_{1}^{\stackrel{\rightharpoonup}{l}}=l!_{\mathbbm{q}^{-s}q^{-s}}p_{0}^{l},\;\;\;\;
X_{1}^{\stackrel{\leftharpoonup}{l}}=\mathbbm{q}^{sl'(l-l')/n}l!_{\mathbbm{q}^{-s}q^{-s}}p_{0}^{l}.$$
\end{lemma}

\begin{proof} By induction on $l.$ Firstly, we have
\begin{eqnarray*}
X_{1}\cdot X_{1}&=&[g . X_{1}][X_{1} . \unit]+[X_{1} .
g][\unit . X_{1}]\\
&=&(1+\mathbbm{q}^{-s}q^{-s})X_{2}X_{1}
=2_{\mathbbm{q}^{-s}q^{-s}}X_{2}X_{1}.
\end{eqnarray*}
Assume $l=an+i$ with $0 \le i \le n-1$ and
$X_{1}^{\stackrel{\rightharpoonup}{l-1}}=(l-1)!_{\mathbbm{q}^{-s}q^{-s}}p_{0}^{l-1},$
then
\begin{eqnarray*}
X_{1}^{\stackrel{\rightharpoonup}{l}}&=&X_{1}^{\stackrel{\rightharpoonup}{l-1}}
\cdot X_{1}
=(l-1)!_{\mathbbm{q}^{-s}q^{-s}}p_{0}^{l-1}\cdot X_{1}\\
&=&(l-1)!_{\mathbbm{q}^{-s}q^{-s}}([g^{i} . X_{1}][X_{i} .
\unit]\cdots[X_{1} . \unit]\\
&&+\cdots+[X_{i} . g]\cdots [X_{1} .
g][\unit . X_{1}] )\\
&=&(l-1)!_{\mathbbm{q}^{-s}q^{-s}}(1+\mathbbm{q}^{-s}q^{-s}+\cdots +
\mathbbm{q}^{-(l-1)s}q^{-(l-1)s})p_{0}^{l}\\
&=& l!_{\mathbbm{q}^{-s}q^{-s}}p_{0}^{l}.
\end{eqnarray*}

Similarly, for $l=an+i$ as above we have

\begin{eqnarray*}
X_{1}^{\stackrel{\leftharpoonup}{l+1}}&=&X_{1}\cdot
X_{1}^{\stackrel{\leftharpoonup}{l}}
=\mathbbm{q}^{sia}l!_{\mathbbm{q}^{-s}q^{-s}} X_{1}\cdot p_{0}^{l}\\
&=&\mathbbm{q}^{sia}l!_{\mathbbm{q}^{-s}q^{-s}}((\mathbbm{q}^{-s}q^{-s})^{i}+\cdots
(\mathbbm{q}^{-s}q^{-s})^{1}+1\\
&&+\mathbbm{q}^{s}n_{\mathbbm{q}^{-s}q^{-s}}
+\cdots +\mathbbm{q}^{sa}n_{\mathbbm{q}^{-s}q^{-s}})p_{0}^{l+1}\\
&=&\mathbbm{q}^{sia}l!_{\mathbbm{q}^{-s}q^{-s}}\mathbbm{q}^{sa}(l+1)_{\mathbbm{q}^{-s}q^{-s}}p_{0}^{l+1}\\
&=& \mathbbm{q}^{s(i+1)a}(l+1)!_{\mathbbm{q}^{-s}q^{-s}}p_{0}^{l+1}.
\end{eqnarray*}
\end{proof}

\begin{lemma} For all $0\leq i,j\leq n-1$ and all non-negative integers $a,b$,
$$X_{1}^{\stackrel{\rightharpoonup}{an+i}} \cdot X_{1}^{\stackrel{\leftharpoonup}{bn+j}}=\left \{
\begin{array}{ll} \mathbbm{q}^{-(a+1)(i+j)s} X_{1}^{\stackrel{\leftharpoonup}{(a+b)n+(i+j)}} & \;\;\;\;
\emph{if} \; i+j> n-1\\ \mathbbm{q}^{-a(i+j)s}
X_{1}^{\stackrel{\leftharpoonup}{(a+b)n+(i+j)}} & \;\;\;\;\emph{if}
\; i+j\leq n-1.
\end{array}\right. $$
\end{lemma}
\begin{proof} Consider the multiplication of $X_{1}^{\stackrel{\rightharpoonup}{an+i-1}},$ $X_1$ and
$X_{1}^{\stackrel{\leftharpoonup}{bn+j}}.$ By the
quasi-associativity axiom (2.1) of Majid algebras, we have the
following equality
$$X_{1}^{\stackrel{\rightharpoonup}{an+i}} \cdot X_{1}^{\stackrel{\leftharpoonup}{bn+j}}
=\frac{\Phi(\unit,\unit,\unit)}{\Phi(g^{i-1},g,g^{j})}X_{1}^{\stackrel{\rightharpoonup}{an+i-1}}
\cdot X_{1}^{\stackrel{\leftharpoonup}{bn+j+1}} \ .$$ By this, it
follows that the coefficient is trivial unless $j=n-1.$ Thus if
$i+j> n-1$, by making use of the previous equation iteratively we
have
\begin{eqnarray*}
X_{1}^{\stackrel{\rightharpoonup}{an+i}} \cdot
X_{1}^{\stackrel{\leftharpoonup}{bn+j}}&=&X_{1}^{\stackrel{\rightharpoonup}{an+(i+j+1-n)}}\cdot
X_{1}^{\stackrel{\leftharpoonup}{bn+n-1}}\\
&=& \mathbbm{q}^{-(i+j)s} X_{1}^{\stackrel{\rightharpoonup}{an+(i+j-n)}}\cdot X_{1}^{\stackrel{\leftharpoonup}{(b+1)n}}\\
&=&\mathbbm{q}^{-(i+j)s} X_{1}^{\stackrel{\rightharpoonup}{(a-1)n+i}}\cdot X_{1}^{\stackrel{\leftharpoonup}{(b+1)n+j}}\\
&=&\cdots\\
&=&\mathbbm{q}^{-a(i+j)s} X_{1}^{\stackrel{\rightharpoonup}{i}}\cdot X_{1}^{\stackrel{\leftharpoonup}{(a+b)n+j}}\\
&=&\mathbbm{q}^{-a(i+j)s} X_{1}^{\stackrel{\rightharpoonup}{i+j+1-n}}\cdot X_{1}^{\stackrel{\leftharpoonup}{(a+b)n+n-1}}\\
&=&\mathbbm{q}^{(-a-1)(i+j)s} X_{1}^{\stackrel{\rightharpoonup}{i+j-n}}\cdot X_{1}^{\stackrel{\leftharpoonup}{(a+b+1)n}}\\
&=&\mathbbm{q}^{-(a+1)(i+j)s}
X_{1}^{\stackrel{\leftharpoonup}{(a+b)n+(i+j)}}.
\end{eqnarray*}

The case of $i+j\leq n-1$ can be proved in the same manner.
\end{proof}

With these preparations, now we can give the product formula.

\begin{proposition} For all non-negative integers $l,m$, we have
$$p_{0}^{l}\cdot p_{0}^{m}=\mathbbm{q}^{sl'(m-m')/n}\binom{l+m}{l}_{\mathbbm{q}^{-s}q^{-s}}p_{0}^{l+m} \ .$$
\end{proposition}

\begin{proof} Assume that $l=an+i,\;m=bn+j$ with $0\leq i,j\leq
n-1$. We only prove the formula when $i+j\leq n-1$ since the case
$i+j > n-1$ can be proved similarly. Indeed, by Lemmas 3.6 and 3.7,
we have
\begin{eqnarray*}
p_{0}^{l}\cdot p_{0}^{m}&=&
\frac{1}{l!_{\mathbbm{q}^{-s}q^{-s}}}\;\frac{1}{\mathbbm{q}^{sjb}m!_{\mathbbm{q}^{-s}q^{-s}}}\;
X_{1}^{\stackrel{\rightharpoonup}{l}} \cdot X_{1}^{\stackrel{\leftharpoonup}{m}}\\
&=&\frac{\mathbbm{q}^{-sjb}\mathbbm{q}^{-a(i+j)s}}{l!_{\mathbbm{q}^{-s}q^{-s}}m!_{\mathbbm{q}^{-s}q^{-s}}}\;
X_{1}^{\stackrel{\leftharpoonup}{l+m}}\\
&=&\frac{\mathbbm{q}^{-sjb}\mathbbm{q}^{-a(i+j)s}}{l!_{\mathbbm{q}^{-s}q^{-s}}m!_{\mathbbm{q}^{-s}q^{-s}}}
\mathbbm{q}^{s(i+j)(a+b)} (l+m)!_{\mathbbm{q}^{-s}q^{-s}}p_{0}^{l+m}\\
&=&\mathbbm{q}^{sib}\binom{l+m}{l}_{\mathbbm{q}^{-s}q^{-s}}p_{0}^{l+m}
\ .
\end{eqnarray*}
\end{proof}

Observe that
\begin{equation} g^i \cdot p_0^l =
\mathbbm{q}^{si(l-l')/n}p_i^l, \quad \quad p_0^l \cdot g^i =
\mathbbm{q}^{-sil}q^{-sil}p_i^l.
\end{equation}
Then the previous product formula can be extended to general paths
as follows. We leave the proof to the interested reader.
\begin{corollary}
For all $0\leq i,j\leq n-1$ and all non-negative integers $l,m$, we
have in $\c Z^n(s,q)$ the following multiplication formula
$$p_i^{l}\cdot p_j^{m}=\mathbbm{q}^{-sjl}q^{-sjl}\mathbbm{q}^{s(i+l')[m+j-(m+j)']/n}
\binom{l+m}{l}_{\mathbbm{q}^{-s}q^{-s}}p_{i+j}^{l+m} \ .$$
\end{corollary}

We remark that the set $\{ \c Z^n(0,q) | q^n=1 \} \cup \{ \c
Z^n(s,q) | 1 \le s \le n-1, \ q \text{ is a primitive root of unity
of order } n^2 \}$ gives an explicit classification of graded Majid
algebras on the path coalgebra $\c Z^n$ by neglecting minor
difference of the quasi-antipodes \cite{d2}. In particular, when
$s=0,$ the set $\{ \c Z^n(0,q) | q^n=1 \}$ gives a classification of
graded Hopf algebras on $\c Z^n,$ which recovers the result in
\cite{c1}; when $s \ne 0,$ the set $\{ \c Z^n(s,q) | 1 \le s \le
n-1, \ q \text{ is a primitive root of unity of order } n^2 \}$
gives a classification of non-trivial graded Majid algebras (i.e.,
not gauge equivalent to Hopf algebras) on $\c Z^n$ by \cite{qha2}.

\subsection{Finite Sub Majid Algebra on $Z^n$} In this subsection we
investigate the possible finite-dimensional graded large sub Majid
algebras of $\c Z^n(s,q).$

\begin{proposition} There is a unique finite-dimensional graded large sub
Majid algebra of $\c Z^n(s,q)$. In addition,  such unique Majid
algebra is generated by $g$ and $X_1.$
\end{proposition}

\begin{proof}
By definition, the smallest graded large sub Majid algebra of $\c
Z^n(s,q)$ is the one generated by the set of vertices and arrows,
which can be given by $g$ and $X_1$ clearly. Assume that the
multiplicative order of $\mathbbm{q}^{-s}q^{-s}$ is $d.$ Then by
Lemmas 3.5 and 3.6, the quasi-algebra generated by $g$ and $X_1$ has
$\{ g^i \cdot X_1^l | 0 \le i \le n-1, 0 \le l \le d-1 \}$ as a
basis. Thus by (3.13) the set $\{p_i^l | 0 \le i \le n-1, 0 \le l
\le d-1 \}$ is also a basis. Clearly the space spanned by these
paths is closed under the coproduct of the path coalgebra, and also
closed under the counit and the quasi-antipode of $\c Z^n(s,q).$
Therefore it is an $nd$-dimensional graded large sub Majid algebra.
We denote this graded Majid algebra by $M(n,s,q).$

It remains to prove that any graded large sub Majid algebra of $\c
Z^n(s,q)$ is infinite-dimensional if it strictly contains
$M(n,s,q).$ Assume that $M$ is such a Majid algebra. Then consider a
nontrivial homogeneous space of degree $l \ge d$ in $M.$ Note that
such a space must be spanned by some paths of length $l.$ Since
there is only one path of length $l$ with fixed source and target,
it follows by the axioms of coalgebras that there are paths of
length $l$ lies in $H.$ Then by (3.13) all the paths of length $l$
lie in $H.$ Consider the coproduct of these paths, it follows that
all the paths of length $d$ must lie in $H.$ Now by making use of
Proposition 3.8 with induction, we have
$(p_0^d)^{\stackrel{\rightharpoonup}{m}} = m! p_0^{md}.$ It follows
that $p_0^{md} \in H$ for all $m \ge 0.$ That means $H$ must be
infinite-dimensional. We are done.
\end{proof}

Now we can conclude that the set $\{ M(n,0,q) | q^n=1 \} \cup \{
M(n,s,q) | 1 \le s \le n-1, \ q \text{ is a primitive root of unity
of order } n^2 \}$ provides a complete classification of
finite-dimensional graded large sub Majid algebra on the Hopf quiver
$Z^n.$ Note that the set $\{ M(n,0,q) | q^n=1 \}$ are the usual Hopf
algebras, called generalized Taft algebras in \cite{hcz}. This gives
a classification of finite-dimensional graded Hopf algebras on the
quiver $Z^n.$ While the set $\{ M(n,s,q) | 1 \le s \le n-1, \ q
\text{ is a primitive root of unity of order } n^2 \}$ gives a
classification of finite-dimensional graded non-trivial Majid
algebras on $Z^n.$

It is also worthy to remark that the underlying coalgebras of these
Majid algebras are truncated subcoalgebras of $\c Z^n,$ namely they
have a basis consisting of all the paths of length smaller than some
fixed integer $d \ge 2.$ In particular, they are monomial in the
sense of \cite{chyz}. The set $\{ M(n,0,q) | q^n=1 \} \cup \{
M(n,s,q) | 1 \le s \le n-1, \ q \text{ is a primitive root of unity
of order } n^2 \}$ also gives a classification of connected monomial
graded Majid algebras, which contains the classification result of
monomial Hopf algebras in \cite{chyz}.

\begin{corollary}
Let $\c Z^n(d)$ denote the truncated sub coalgebra of $\c Z^n$
spanned by the paths of length smaller than some fixed integer $d
\ge 2.$ Then $\c Z^n(d)$ admits a graded Majid algebra structure if
and only if $d | n,$ or $d=\frac{n^2}{(s,n^2)}$ for some $1 \le s
\le n-1.$ In this case, the Majid algebra on $\c Z^n(d)$ is gauge
equivalent to a Hopf algebra if and only if $d | n.$
\end{corollary}

\subsection{Classification Results}

Now we are ready to give the main result. In the following we always
assume that, for a pointed Majid algebra $M$ with set of group-like
elements $G(M),$ its quasi-antipode $(\mathcal{S},\alpha,\beta)$
satisfies $\S(g)=g^{-1}, \ \alpha(g)=1, \
\beta(g)=1/\Phi(g,g^{-1},g)$ for all $g \in G(M).$ The observation
of Drinfeld \cite{d2} guarantees that this assumption is harmless.

\begin{theorem} Suppose that $M \ne \c$ is a finite-dimensional graded pointed
Majid algebra with connected underlying coalgebra. If $M$ is of
finite corepresentation type, then $M \cong M(n,s,q)$ for some
positive integers $n \ge 2$ and $0 \le s \le n-1,$ and $q$ is an
$n$-th root of unity if $s=0,$ or some primitive root of unity of
order $n^2$ if $s \ne 0.$
\end{theorem}

\begin{proof}
By the assumption, the corresponding Hopf quiver of $M$ is $Z^n$ for
some $n \ge 2$ according to Corollary 3.2. Now by the Gabriel type
theorem for Majid algebra and Subsection 3.4, $M$ can be viewed as a
graded large sub Majid algebra of some $\c Z^n(s,q).$ By Proposition
3.10, there is only one possible finite-dimensional graded large sub
Majid algebra $M(n,s,q).$ So $M$ can only be one of the $M(n,s,q).$

On the other hand, since the $M(n,s,q)$ are truncated sub coalgebra
of $\c \Z^n,$ so they are of finite corepresentation type as
mentioned in Lemma 3.1. This completes the proof.
\end{proof}

With a help of this theorem, in the following we deal with not
necessarily graded situation by making use of deformation theory
\cite{ger} and a deep result from geometric methods of
representation theory \cite{gab3}. Without loss of generality, let
$M$ be a finite-dimensional pointed Majid algebra with connected
underlying coalgebra, that is, its bound quiver is connected. As
before, we exclude the trivial case for $M=\c.$

\begin{corollary}
Keep the above assumption. Then $M$ is of finite corepresentation
type if and only if $\gr M \cong M(n,s,q)$ for some appropriate
$n,s,q.$
\end{corollary}

\begin{proof}
First suppose that $M$ is of finite corepresentation type. Consider
its underlying coalgebra and apply the Gabriel type theorem. Since
$M$ and $\gr M$ share the same quiver, so there is a unique Hopf
quiver $Q$ such that $M$ can be viewed as a large sub coalgebra of
the path coalgebra $\c Q.$ Now Lemma 3.1 applies, that is, when $M$
is of corepresentation type, then the quiver $Q$ can only be $Z^n$
for some integer $n \ge 2.$ According to Proposition 3.10, the
graded version $\gr M$ of $M$ is isomorphic to some $M(n,s,q).$

Conversely, suppose $\gr M \cong M(n,s,q).$ Note that $M$ is a
deformation of $\gr M$ (it is also said that $\gr M$ is a
degeneration of $M$). Then that the Majid algebra $M(n,s,q)$ is of
finite corepresentation type implies that so is $M,$ according to
the famous theorem of Gabriel \cite{gab3} which says that finite
representation type is open.
\end{proof}

To complete the classification of non-graded connected pointed Majid
algebras of finite corepresentation type, it suffices to calculate
all the deformations of $M(n,s,q).$ Note that finite-dimensional
Majid algebras are co-Frobenius (i.e., the dual algebra is
Frobenius) according to \cite{hn,bc}, then by the same argument as
in Section 2 of \cite{chyz} we can conclude that the underlying
coalgebra of a connected pointed Majid algebras of finite
corepresentation type is isomorphic to a truncated sub coalgebra of
$\c Z^n.$ It follows that one only needs to calculate the
coalgebra-preserving deformation of $M(n,s,q),$ which is a quasi
analogue of the preferred deformation of Hopf algebras
\cite{ger-sch}.

\subsection{Some Remarks} We conclude this section with some
remarks.

(1) The preferred deformations for $M(n,0,q)$ were explicitly given
in \cite{chyz}. For $s \ne 0,$ it seems that the preferred
deformations for $M(n,s,q)$ are much more complicated since the
deformation of reassociators gets involved. We leave this problem
for future work.

(2) For not necessarily connected situation, the underlying
coalgebra of a finite-dimensional pointed Majid algebras is a direct
sum of finite copies of some $\c Z^n(d)$ and the Majid algebra is a
crossed product of some deformation of $M(n,s,q)$ with a group
twisted by a three cocycle \cite{qha1}.

(3) Our classification of finite-dimensional pointed Majid algebras
of finite corepresentation type contains the corresponding
classification result for pointed Hopf algebras, which was given in
\cite{ll} by different method.

(4) A standard dualization process gives parallel classification
results for elementary quasi-Hopf algebras of finite representation
type. Some of the dual of $M(n,s,q)$ appeared in previous works of
Etingof and Gelaki \cite{eg1,eg2,eg3,g}.

\section{Tensor categories of finite type}

As an application, we will classify a class of \emph{tensor
categories of finite type} in this section. In addition, some
information of these finite tensor categories are given by making
use of quiver representation theory.

\subsection{Finite Tensor Categories}

By a tensor category we mean an abelian rigid monoidal category over
$\c$ in which the neutral object $\mathbf{1}$ is simple. A tensor
category $\mathscr{C}$ is said to be finite if

\begin{enumerate}
  \item $\mathscr{C}$ has finitely many simple objects,
  \item any object has finite length, and
  \item any simple object admits a projective cover.
\end{enumerate}

For a finite tensor category $\mathscr{C}$, we denote its
Grothendieck ring by $Gr(\mathscr{C}).$ It is a free abelian group
of finite rank, whose basis $S$ is the set of isomorphism classes of
simple objects in $\mathscr{C}.$ Let $X\in \mathscr{C}$, then its
\emph{Frobenius-Perron dimension} $d_{+}(X)$ is defined to be the
largest non-negative real eigenvalue of the matrix of left
multiplication in the Grothendieck ring by $X$ under the basis $S.$
For more knowledge about finite tensor categories, see \cite{eno,eo,
ce} and references therein.

The following result of Etingof and Ostrik \cite{eo} provides the
close relation between finite tensor categories and finite
quasi-quantum groups.

\begin{lemma}\textbf{\emph{(Etingof-Ostrik)}} For a finite tensor category $\mathscr{C}$,
it is tensor equivalent to the representation category $\Rep H$ of a
finite-dimensional quasi-Hopf algebra $H$ if and only if the
Frobenius-Perron dimensions of objects in $\mathscr{C}$ are
integers.
\end{lemma}

Since Majid algebras are dual of quasi-Hopf algebras, so a finite
tensor category $\mathscr{C}$ is tensor equivalent to the
corepresentation category $\Corep M$ of a finite-dimensional Majid
algebra $M$ if and only if the Frobenius-Perron dimensions of
objects in $\mathscr{C}$ are integers.

\subsection{Tensor Categories of Finite Type}

The simplest finite tensor categories are of course the semisimple
ones. After the semisimple situation, the simplest ones are those
having finitely many isomorphism classes of indecomposable objects
in view of the Krull-Remak-Schmidt property of finite tensor
categories. A finite tensor category $\mathscr{C}$ with this
property must be equivalent to the representation category $\Rep A$
of a finite-dimensional algebra $A$ of finite representation type,
or the corepresentation category $\Corep C$ of a finite-dimensional
coalgebra $C$ of finite corepresentation type. Inspired by this,
such tensor categories are said to be of \emph{finite type}.

By Lemma 4.1, the classification of tensor categories of finite type
whose objects have integer Frobenius-Perron dimensions is equivalent
to the classification of finite-dimensional quasi-Hopf algebras of
finite representation type, or finite-dimensional Majid algebras of
finite corepresentation type. Now the results of Section 3 can be
applied to classify some class of such tensor categories.

\subsection{Some Classification Results}

Let $\mathscr{C}$ be a finite tensor category and assume that
$\mathscr{C} = \Corep C$ as an abelian category for a
finite-dimensional coalgebra $C.$ We say that $\mathscr{C}$ is
connected if the coalgebra $C$ is connected, that is, the dual
algebra $C^*$ of $C$ is indecomposable \cite{ars}. To avoid the
trivial case, in the following we always assume that $\mathscr{C}$
has at least two simple objects.

\begin{theorem}
Assume that $\mathscr{C}$ is a tensor category of finite type. If
every simple object of $\mathscr{C}$ has Frobenius-Perron dimension
1, then as an abelian category $\mathscr{C}$ is equivalent to direct
product of finite copies of $\Corep\c Z^n(d)$ with $d | n,$ or
$d=\frac{n^2}{(s,n^2)}$ for some $1 \le s \le n-1.$
\end{theorem}

\begin{proof}
By the assumption, first of all we have $\mathscr{C} = \Corep C$ as
an abelian category for some finite-dimensional pointed coalgebra
$C.$ By the property of Frobenius-Perron dimension, the fact that
every simple object of $\mathscr{C}$ has Frobenius-Perron dimension
1 implies all objects of $\mathscr{C}$ have integer Frobenius-Perron
dimensions. Then $C$ must be a pointed Majid algebra of finite
corepresentation type by Lemma 4.1. Now the theorem follows
immediately from Subsections 3.6 and 3.7.
\end{proof}

The corepresentation category  $\Corep\c Z^n(d)$ is well understood,
see for example \cite{ass,ars} and also \cite{hcz} for its
Auslander-Reiten quiver. In particular, $\c Z^n(d)$ is Nakayama
\cite{ass,ars}. Therefore, we have

\begin{corollary}
Assume that $\mathscr{C}$ is a tensor category of finite type and
that every simple object of $\mathscr{C}$ has Frobenius-Perron
dimension 1, then every indecomposable object $X \in \mathscr{C}$ is
uniserial, that is, the set of sub objects of $X$ is totally ordered
by inclusion.
\end{corollary}

Combine Corollary 3.13 and Lemma 4.1 we have

\begin{corollary}
Assume that $\mathscr{C}$ is a connected tensor category of finite
type. If every simple object of $\mathscr{C}$ has Frobenius-Perron
dimension 1, then as a tensor category $\mathscr{C}$ is equivalent
to $\Corep M$ for some pointed Majid algebras $M$ with $\gr M \cong
M(n,s,q).$
\end{corollary}

\vskip 0.5cm

\noindent{\bf Acknowledgements:} The research of the first author
was supported by the NSF of China under grant number 10601052. The
research of the second author was supported by the NSF of China
under grant number 10801069.

\end{document}